\documentclass[12pt,twoside]{article}
\usepackage[hmargin=0.8in,vmargin=0.9in]{geometry}
\usepackage{fancyhdr}
\usepackage{graphicx}
\usepackage{subfigure}
\usepackage{amssymb}
\usepackage{amsmath}
\usepackage{amsfonts}
\usepackage{theorem}
\usepackage{mathrsfs}
\usepackage{mathtools}
\usepackage{empheq}
\usepackage{bm}
\usepackage{bbm}
\usepackage{color}
\usepackage{setspace}
\usepackage{epstopdf}
\usepackage{float}
\usepackage{picinpar}
\usepackage{extarrows}
\usepackage{multirow}

\usepackage{cite}
\usepackage{nicefrac}
\usepackage{booktabs} % for much better looking tables
\usepackage{array} % for better arrays (eg matrices) in maths
\usepackage{paralist} % very flexible & customisable lists (eg. enumerate/itemize, etc.)
\usepackage{verbatim} % adds environment for commenting out blocks of text & for better verbatim
\usepackage{subfigure} % make it possible to include more than one captioned figure/table in a single float
\usepackage{authblk}
\usepackage[normalem]{ulem}

\newtheorem{thm}{Theorem}[section]

\newtheorem{rmk}[thm]{Remark}

{\theorembodyfont{\rmfamily}}

\newenvironment{DA}{{\flushleft \bf Declarations:}}{}

\newcommand\eref[1]{(\ref{#1})}

\newcommand*\xbar[1]{%
  \hbox{%
    \vbox{%
      \hrule height 0.5pt % The actual bar
      \kern0.4ex%         % Distance between bar and symbol
      \hbox{%
        \kern-0.05em%      % Shortening on the left side
        \ensuremath{#1}%
        \kern-0.00em%      % Shortening on the right side
      }%
    }%
  }%
}

\allowdisplaybreaks[1]

\def\hf {\frac{1}{2}}

\newcommand{\jph}{{j+\frac{1}{2}}}
\newcommand{\jmh}{{j-\frac{1}{2}}}
\newcommand{\kph}{{k+\frac{1}{2}}}
\newcommand{\kmh}{{k-\frac{1}{2}}}

\newcommand{\dx}{\Delta x}
\newcommand{\dy}{\Delta y}
\newcommand{\dt}{\Delta t}

\numberwithin{equation}{section}
\numberwithin{figure}{section}
\numberwithin{table}{section}

\graphicspath{{figures/}{./}}

\title{New Low-Dissipation Central-Upwind Schemes. Part II}
\author{Shaoshuai Chu\thanks{Department of Mathematics and Shenzhen International Center for Mathematics, Southern University of Science and
Technology, Shenzhen, 518055, China; {\tt chuss@mail.sustech.edu.cn}}, Alexander Kurganov\thanks{Department of Mathematics, Shenzhen
International Center for Mathematics, and Guangdong Provincial Key Laboratory of Computational Science and Material Design, Southern
University of Science and Technology, Shenzhen, 518055, China; {\tt alexander@sustech.edu.cn}}, and Ruixiao Xin\thanks{Department of
Mathematics, Southern University of Science and Technology, Shenzhen, 518055, China; {\tt 12331009@mail.sustech.edu.cn}}}

\begin{document}

\date{}
\maketitle
\begin{abstract}
The low-dissipation central-upwind (LDCU) schemes have been recently introduced in  [{\sc A. Kurganov and R. Xin, J. Sci. Comput., 96
(2023), Paper No. 56}] as a modification of the central-upwind (CU) schemes from [{\sc A. Kurganov and C. T. Lin, Commun. Comput. Phys., 2
(2007), pp. 141-163}]. The LDCU schemes achieve much higher resolution of contact waves and many (two-dimensional) structures resulting from
complicated wave interaction. However, the LDCU schemes sometimes produce more oscillatory results compared with the CU schemes, especially
near the computational domain boundaries.

In this paper, we propose a very simple---yet systematic---modification of the LDCU schemes, which completely eliminates the aforementioned
oscillations almost without affecting the quality of the computed solution.
\end{abstract}



\section{Introduction}
We consider hyperbolic systems of conservation laws, which in the one-dimensional (1-D) case read as
\begin{equation}
\bm U_t+\bm F(\bm U)_x=\bm0,
\label{1.1}
\end{equation}
and in the two-dimensional (2-D) case read as
\begin{equation}
\bm U_t+\bm F(\bm U)_x+\bm G(\bm U)_y=\bm0,
\label{1.2}
\end{equation}
where $x$ and $y$ are spatial variables, $t$ is the time, $\bm U\in\mathbb R^d$ is a vector of unknowns, and $\bm F$ and $\bm G$ are the
$x$- and $y$-directional fluxes, respectively.

We focus on finite-volume central-upwind (CU) schemes, which were introduced in \cite{Kurganov01,Kurganov00,KT2002} as a ``black-box''
solver for general hyperbolic systems \eref{1.1} and \eref{1.2}. Even though the CU schemes are quite accurate, efficient, and robust for a
wide variety of hyperbolic systems, higher resolution of the numerical solutions can be achieved by further reducing numerical dissipation
present in the CU schemes; see, e.g., \cite{Kurganov07,CCHKL_22,KX_22}. The low-dissipation CU (LDCU) schemes have been recently proposed in
\cite{KX_22} as a modification of the CU scheme from \cite{Kurganov07}. The LDCU schemes achieve much higher resolution of contact waves and
many (two-dimensional) structures resulting from complicated wave interaction. However, the LDCU schemes sometimes produce more oscillatory
results compared with the CU schemes, especially near the computational domain boundaries.

In order to suppress these spurious oscillations, we modify both 1-D and 2-D LDCU schemes. The proposed modifications are very simple, yet
systematic as they are based on a more accurate projection of the evolved solution onto the original (uniform) mesh. The new LDCU schemes
are developed for the 1-D and 2-D Euler equations of gas dynamics and tested on several numerical examples. The obtained numerical results
demonstrate that the new schemes contain almost the same small amount of numerical dissipation as the LDCU schemes from \cite{KX_22} but
produce substantially ``cleaner'', non-oscillatory computed solutions.

\section{New LDCU Schemes}\label{sec2}
We follow the derivation of the LDCU scheme in \cite{KX_22} and use precisely the same notation as in \cite{KX_22}.

We cover the computational domain with the finite volume cells $C_j=[x_\jmh,x_\jph]$, $j=1,\ldots,N$, which are assumed to be uniform, that
is, $x_\jph-x_\jmh\equiv\dx$. We then assume that the solution, realized in terms of its cell averages $\,\xbar{\bm U}^{\,n}_j$, is
available at a certain time level $t=t^n$ and reconstruct a second-order piecewise linear interpolant
$\sum_{j}[\xbar{\bm U}^{\,n}_j+(\bm U_x)^n_j(x-x_j)]{\cal X}_{C_j}$, where $\cal X$ denotes the characteristic function of the corresponding
intervals and $(\bm U_x)^n_j$ are the slopes obtained using a nonlinear limiter. In the numerical experiments reported in \S\ref{sec3}, we
have used a generalized minmod limiter with the minmod parameter $\theta=1.3$; see, e.g., \cite{LN,NT,Swe,vLeV}.

We then evaluate the local speeds of propagation $a^\pm_\jph$, introduce the corresponding points $x^n_{\jph,\ell}:=x_\jph+a^-_\jph\dt^n$
and $x^n_{\jph,r}:=x_\jph+a^+_\jph\dt^n$, and integrate the system \eref{1.1} over the space-time control volumes consisting of the
``smooth'' $[x_{\jmh,r},x_{\jph,\ell}]\times[t^n,t^{n+1}]$ and ``nonsmooth'' $[x_{\jph,\ell},x_{\jph,r}]\times[t^n,t^{n+1}]$ areas, where
$t^{n+1}:=t^n+\dt^n$. This way the solution is evolved in time and upon the completion of the evolution step, we obtain the intermediate
cell averages $\xbar{\bm U}^{\rm\,int}_\jph$ (see \cite[Eq. (2.7)]{KX_22}) and $\xbar{\bm U}^{\rm\,int}_j$ (see \cite[Eq. (2.10)]{KX_22}).

Next, the intermediate solution, which is realized in terms of $\{\xbar{\bm U}^{\,\rm int}_j\}$ and $\{\xbar{\bm U}^{\,\rm int}_\jph\}$, is
projected onto the original grid. To this end, we need to construct the interpolant
\begin{equation}
\widetilde{\bm U}^{\,\rm int}(x)=\sum_j\left\{\widetilde{\bm U}^{\,\rm int}_\jph(x){\cal X}_{[x_{\jph,\ell},x_{\jph,r}]}+
\xbar{\bm U}^{\,\rm int}_j{\cal X}_{[x_{\jmh,r},x_{\jph,\ell}]}\right\}.
\label{2.2.4}
\end{equation}
In order to develop the original LDCU scheme in \cite{KX_22}, we have used
\begin{equation}
\widetilde{\bm U}^{\,\rm int}_\jph(x)=\left\{\begin{aligned}
&\xbar{\bm U}^{\,\rm int,L}_\jph,&&x<x_\jph,\\
&\xbar{\bm U}^{\,\rm int,R}_\jph,&&x>x_\jph,
\end{aligned}\right.
\label{2.2.5}
\end{equation}
where $\widetilde{\bm U}^{\,\rm int}_\jph(x)$ is discontinuous at $x=x_\jph$. The values $\,\xbar{\bm U}^{\,\rm int,L}_\jph$ and
$\,\xbar{\bm U}^{\,\rm int,R}_\jph$ are determined based on the conservation requirement
\begin{equation}
a^+_\jph\,\xbar{\bm U}^{\,\rm int,R}_\jph-a^-_\jph\,\xbar{\bm U}^{\,\rm int,L}_\jph=(a^+_\jph-a^-_\jph)\,\xbar{\bm U}^{\,\rm int}_\jph.
\label{3.0a}
\end{equation}
and $d$ additional degrees of freedom, which can be used to design an accurate projection. The way these degrees freedom are utilized
depends on the problem at hand. One, however, can introduce an additional degree of freedom, which may be used to further improve the
accuracy of the projection step.

\subsection{Modification of the Projection Step}
\begin{figure}[ht!]
\centerline{\includegraphics[width=9.5cm]{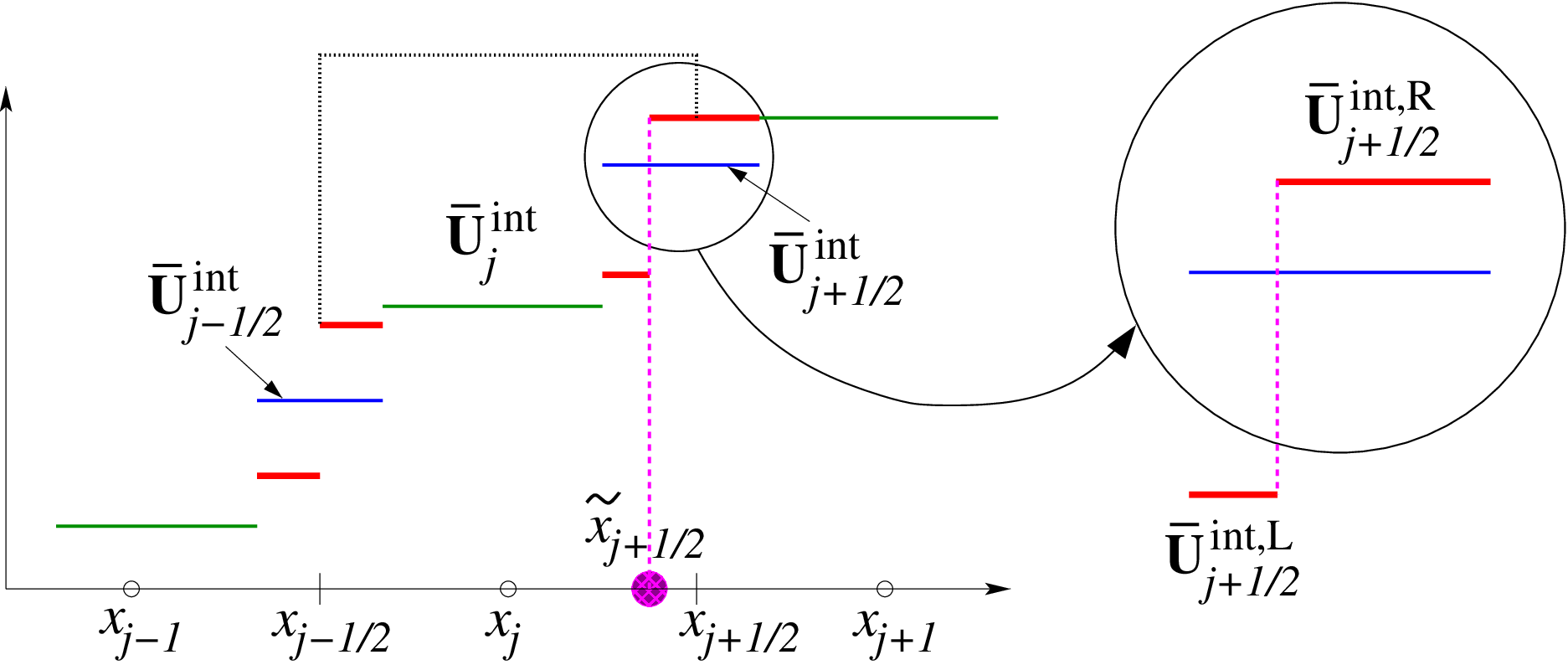}}
\caption{\sf The new projection step.\label{fig3}}
\end{figure}

We propose to replace \eref{2.2.5} with
\begin{equation}
\widetilde{\bm U}^{\,\rm int}_\jph(x)=\left\{\begin{aligned}
&\xbar{\bm U}^{\,\rm int,L}_\jph,&&x<\widetilde x_\jph,\\
&\xbar{\bm U}^{\,\rm int,R}_\jph,&&x>\widetilde x_\jph,
\end{aligned}\right.
\label{2.2.6}
\end{equation}
where $\widetilde x_\jph=x_\jph+\widetilde u_\jph\dt^n$ and $\widetilde u_\jph\in(a^-_\jph,a^+_\jph)$ represents an additional degree of
freedom. The conservation requirements now give
\begin{equation}
(a^+_\jph-\widetilde u_\jph)\,\xbar{\bm U}^{\,\rm int,R}_\jph+(\widetilde u_\jph-a^-_\jph)\,\xbar{\bm U}^{\,\rm int,L}_\jph=
(a^+_\jph-a^-_\jph)\,\xbar{\bm U}^{\,\rm int}_\jph;
\label{2.2.7}
\end{equation}
compare with \eref{3.0a}, which is obtained from \eref{2.2.7} if $\widetilde u_\jph=0$. The modified projection step is outlined in Figure
\ref{fig3}; compare it with \cite[Fig. 3]{KX_22}.

In order to complete the derivation of the new LDCU scheme, we need to consider a particular hyperbolic system. As in \cite{KX_22}, we study
the Euler equation of gas dynamics.

\subsection{New LDCU Scheme for the 1-D Euler Equations of Gas Dynamics}\label{sec22}
The 1-D Euler equations of gas dynamics reads as \eref{1.1} with $\bm U=(\rho,\rho u,E)^\top$ and
$\bm F(\bm U)=(\rho u,\rho u^2+p,u(E+p))^\top$, where $\rho$ is the density, $u$ is the velocity, $p$ is the pressure, and $E$ is the total
energy. The system is closed using the equation of states (EOS), which is in the case of ideal gas is $p=(\gamma-1)\big[E-\hf\rho u^2\big]$,
$\gamma={\rm Const}$.

In order to derive the new LDCU scheme for the 1-D Euler equations, we first follow \cite[\S3.1]{KX_22} to obtain the point values
$\rho^\pm_\jph$, $(\rho u)^\pm_\jph$, and $E^\pm_\jph$, evaluate the corresponding point values $u^\pm_\jph$ and $p^\pm_\jph$, and estimate
the one-sided local speeds of propagation $a^\pm_\jph$.

We then proceed with the evolution of the subcell averages $\,\xbar{\bm U}^{\,\rm int,L}_\jph=
\big(\,\xbar\rho^{\,\rm int,L}_\jph,(\xbar{\rho u})^{\,\rm int,L}_\jph,\xbar E^{\,\rm int,L}_\jph\big)^\top$ and
$\xbar{\bm U}^{\,\rm int,R}_\jph=
\big(\,\xbar\rho^{\,\rm int,R}_\jph,(\xbar{\rho u})^{\,\rm int,R}_\jph,\xbar E^{\,\rm int,R}_\jph\big)^\top$ required in \eref{2.2.6}.
Similarly to \cite[\S3]{KX_22}, we enforce the continuity of $u$ and $p$ across the cell interfaces by setting
\begin{equation}
\begin{aligned}
u^{\rm int,L}_\jph:=\frac{(\xbar{\rho u})^{\,\rm int,L}_\jph}{\xbar\rho^{\,\rm int,L}_\jph}=
\frac{(\xbar{\rho u})^{\,\rm int,R}_\jph}{\xbar\rho^{\,\rm int,R}_\jph}=:u^{\rm int,R}_\jph,\quad
\xbar E^{\,\rm int,L}_\jph-\frac{\big((\xbar{\rho u})^{\,\rm int,L}_\jph\big)^2}{2\,\xbar\rho^{\,\rm int,L}_\jph}=
\xbar E^{\,\rm int,R}_\jph-\frac{\big((\xbar{\rho u})^{\,\rm int,R}_\jph\big)^2}{2\,\xbar\rho^{\,\rm int,R}_\jph}.
\end{aligned}
\label{3.5}
\end{equation}
Next, \eref{3.5} together with the conservation requirement \eref{2.2.7} applied to $\rho$, $\rho u$, and $E$,
\begin{equation}
\begin{aligned}
(a^+_\jph-\widetilde u_\jph)\,\xbar\rho^{\,\rm int,R}_\jph+(\widetilde u_\jph-a^-_\jph)\,\xbar\rho^{\,\rm int,L}_\jph&
=\big(a^+_\jph-a^-_\jph\big)\,\xbar\rho^{\,\rm int}_\jph,\\
(a^+_\jph-\widetilde u_\jph)(\xbar{\rho u})^{\,\rm int,R}_\jph+(\widetilde u_\jph-a^-_\jph)(\xbar{\rho u})^{\,\rm int,L}_\jph&
=\big(a^+_\jph-a^-_\jph\big)(\xbar{\rho u})^{\,\rm int}_\jph,\\
(a^+_\jph-\widetilde u_\jph)\,\xbar E^{\,\rm int,R}_\jph+(\widetilde u_\jph-a^-_\jph)\,\xbar E^{\,\rm int,L}_\jph&
=\big(a^+_\jph-a^-_\jph\big)\,\xbar E^{\,\rm int}_\jph,
\end{aligned}
\label{3.5a}
\end{equation}
yield a system of five algebraic equations \eref{3.5}--\eref{3.5a}, which we solve for $(\xbar{\rho u})^{\,\rm int,L}_\jph$,
$(\xbar{\rho u})^{\,\rm int,R}_\jph$, $\xbar E^{\,\rm int,L}_\jph$, and $\xbar E^{\,\rm int,R}_\jph$, and express these quantities in terms
of $\,\xbar\rho^{\,\rm int,L}_\jph$ and $\,\xbar\rho^{\,\rm int,R}_\jph$:
\begin{equation}
\begin{aligned}
&(\xbar{\rho u})^{\,\rm int,L}_\jph=\,\xbar\rho^{\,\rm int,L}_\jph u^{\rm int}_\jph,~&&
(\xbar{\rho u})^{\,\rm int,R}_\jph=\,\xbar\rho^{\,\rm int,R}_\jph u^{\rm int}_\jph,\\
&\xbar E^{\,\rm int,L}_\jph=\,\xbar E^{\,\rm int}_\jph+\frac{\xbar\rho^{\,\rm int,L}_\jph-\rho^{\rm int}_\jph}{2}
\big(u^{\rm int}_\jph\big)^2,~&&
\xbar E^{\,\rm int,R}_\jph=\,\xbar E^{\,\rm int}_\jph+\frac{\xbar{\rho}^{\,\rm int,R}_\jph-\rho^{\rm int}_\jph}{2}
\big(u^{\rm int}_\jph\big)^2,
\end{aligned}
\label{3.6}
\end{equation}
where $u^{\rm int}_\jph:=(\xbar{\rho u})^{\,\rm int}_\jph/\,\xbar\rho^{\,\rm int}_\jph$.

Note that the first two equations in \eref{3.6} give $u^{\rm int,L}_\jph=u^{\rm int}_\jph=u^{\rm int,R}_\jph$, which suggests that in this
piecewise constant solution approximation, the velocity is constant in $[x_{\jph,\ell},x_{\jph,r}]$. We therefore use the same velocity
value and set $\widetilde u_\jph=u^{\rm int}_\jph$.

Finally, we follow the steps in \cite[\S3]{KX_22}, where we make the difference
$\,\xbar\rho^{\,\rm int,R}_\jph-\,\xbar\rho^{\,\rm int,L}_\jph$ as large as possible without creating any new local extrema. To this end, we
denote by
\begin{equation*}
S^-_\jph:=(u^{\rm int}_\jph-a^-_\jph)\big(\,\xbar\rho^{\,\rm int}_\jph-\rho^{\rm int}_{\jph,\ell}\big)\quad\mbox{and}\quad
S^+_\jph:=(a^+_\jph-u^{\rm int}_\jph)\big(\rho^{\rm int}_{\jph,r}-\,\xbar\rho^{\,\rm int}_\jph\big),
\end{equation*}
where the point values $\rho^{\rm int}_{\jph,\ell}$ and $\rho^{\rm int}_{\jph,\ell}$ were introduced in \cite[Eq. (2.14)]{KX_22}, and then
determine $\,\xbar\rho^{\,\rm int,L}_\jph$ and $\,\xbar\rho^{\,\rm int,R}_\jph$ by
\begin{equation}
\begin{aligned}
&\xbar\rho^{\,\rm int,L}_\jph=\,\xbar\rho^{\,\rm int}_\jph+\frac{\delta_\jph}{a^-_\jph-u^{\rm int}_\jph},\\
&\xbar\rho^{\,\rm int,R}_\jph=\,\xbar\rho^{\,\rm int}_\jph+\frac{\delta_\jph}{a^+_\jph-u^{\rm int}_\jph},
\end{aligned}
\qquad\delta_\jph:={\rm minmod}\big(S^-_\jph,S^+_\jph\big);
\label{3.14}
\end{equation}
compare these formulae with \cite[Eq. (3.14)]{KX_22}. Here, ${\rm minmod}(a,b):=\frac{1}{2}({\rm sgn}(a)+{\rm sgn}(b))\min(|a|,|b|)$.

We then substitute \eref{3.14} into \eref{3.6} and obtain
\begin{equation}
\begin{aligned}
&(\xbar{\rho u})^{\,\rm int,L}_\jph=(\xbar{\rho u})^{\,\rm int}_\jph+\frac{\delta_\jph}{a^-_\jph-u^{\rm int}_\jph}u^{\rm int}_\jph,~&&
(\xbar{\rho u})^{\,\rm int,R}_\jph=(\xbar{\rho u})^{\,\rm int}_\jph+\frac{\delta_\jph}{a^+_\jph-u^{\rm int}_\jph}u^{\rm int}_\jph,\\
&\xbar E^{\,\rm int,L}_\jph=\,\xbar E^{\,\rm int}_\jph+\frac{\delta_\jph}{2(a^-_\jph-u^{\rm int}_\jph)}\big(u^{\rm int}_\jph\big)^2,~&&
\xbar E^{\,\rm int,R}_\jph=\,\xbar E^{\,\rm int}_\jph+\frac{\delta_\jph}{2(a^+_\jph-u^{\rm int}_\jph)}\big(u^{\rm int}_\jph\big)^2.
\end{aligned}
\label{3.10f}
\end{equation}

\subsubsection{Fully Discrete Scheme}
We now derive a new fully discrete LDCU scheme based on the new projection step. To this end, we integrate the piecewise constant
interpolant \eref{2.2.4}, \eref{2.2.6} over the cell $C_j$ and obtain
\begin{equation}
\begin{aligned}
\xbar{\bm U}^{\,n+1}_j&=\,\xbar{\bm U}^{\,\rm int}_j+\frac{\dt^n}{\dx}\Big[a^+_\jmh\big(\,\xbar{\bm U}^{\,\rm int,R}_\jmh-
\,\xbar{\bm U}^{\,\rm int}_j\big)-a^-_\jph\big(\,\xbar{\bm U}^{\,\rm int,L}_\jph-\,\xbar{\bm U}^{\,\rm int}_j\big)\\
&+\max(u^{\rm int}_\jmh,0)(\,\xbar{\bm U}^{\,\rm int,L}_\jmh-\,\xbar{\bm U}^{\,\rm int,R}_\jmh)-
\max(u^{\rm int}_\jph,0)(\,\xbar{\bm U}^{\,\rm int,L}_\jph-\,\xbar{\bm U}^{\,\rm int,R}_\jph)\Big]\\
&\hspace*{-1.25cm}\stackrel{\eref{3.14},\eref{3.10f}}{=}\,\xbar{\bm U}^{\rm\,int}_j+\frac{\dt^n}{\dx}
\Big[a^+_\jmh\big(\,\xbar{\bm U}^{\,\rm int}_\jmh-\xbar{\bm U}^{\,\rm int}_j\big)
-a^-_\jph\big(\xbar{\bm U}^{\,\rm int}_\jph-\xbar{\bm U}^{\,\rm int}_j\big)
+\alpha^{\rm int}_\jmh{\bm\delta}_\jmh-\alpha^{\rm int}_\jph{\bm\delta}_\jph\Big],
\end{aligned}
\label{3.15}
\end{equation}
where
\begin{equation}
\alpha^{\rm int}_\jph=\left\{\begin{aligned}
&\frac{a^+_\jph}{a^+_\jph-u^{\rm int}_\jph}&&\mbox{if}~u^{\rm int}_\jph<0,\\
&\frac{a^-_\jph}{a^-_\jph-u^{\rm int}_\jph}&&\mbox{otherwise},
\end{aligned}\right.
%\qquad\bm\delta_\jph=\delta_\jph\Bigg(1,u^{\rm int}_\jph,\frac{\big(u^{\rm int}_\jph\big)^2}{2}\,\Bigg)^\top.
\qquad\bm\delta_\jph=\delta_\jph\begin{pmatrix}1\\u^{\rm int}_\jph\\[0.8ex]\dfrac{1}{2}\big(u^{\rm int}_\jph\big)^2\end{pmatrix}.
\label{3.15a}
\end{equation}

\subsubsection{Semi-Discrete Scheme}\label{sec222}
We now pass to the semi-discrete limit $\dt^n\to0$ in \eref{3.15}--\eref{3.15a} and proceed as in \cite[\S3.1.2]{KX_22} to end up with the
following semi-discretization
\begin{equation}
\frac{\rm d}{\rm dt}\,\xbar{\bm U}_j(t)=-\frac{\bm{{\cal F}}_\jph(t)-{\bm{{\cal F}}}_\jmh(t)}{\dx},
\label{2.17}
\end{equation}
where ${\bm{{\cal F}}}_\jph$ are the modified LDCU numerical fluxes given by
\begin{equation}
\bm{{\cal F}}_\jph=\frac{a^+_\jph\bm F\big(\bm U^-_\jph\big)-a^-_\jph\bm F\big(\bm U^+_\jph\big)}{a^+_\jph-a^-_\jph}+
\frac{a^+_\jph a^-_\jph}{a^+_\jph-a^-_\jph}\left(\bm U^+_\jph-\bm U^-_\jph\right)+\bm q_\jph,
\label{2.18}
\end{equation}
and
\begin{equation}
%\bm q_\jph=\alpha^*_\jph q^\rho_\jph\left(1,u^*_\jph,\frac{(u^*_\jph)^2}{2}\right)^\top,
\bm q_\jph=\alpha^*_\jph q^\rho_\jph\begin{pmatrix}1\\u^*_\jph\\[0.8ex]\dfrac{1}{2}\big(u^*_\jph\big)^2\end{pmatrix}
\label{2.19}
\end{equation}
is a modified ``built-in'' anti-diffusion term with
\begin{equation}
\begin{aligned}
&\rho^*_\jph=\frac{a^+_\jph\rho^+_\jph-a^-_\jph\rho^-_\jph-\left[(\rho u)^+_\jph-(\rho u)^-_\jph\right]}{a^+_\jph-a^-_\jph},\\
&(\rho u)^*_\jph=\frac{a^+_\jph(\rho u)^+_\jph-a^-_\jph(\rho u)^-_\jph-\left[\rho^+_\jph\big(u^+_\jph\big)^2+p^+_\jph-\rho^-_\jph
\big(u^-_\jph\big)^2-p^-_\jph\right]}{a^+_\jph-a^-_\jph},\\
&u^*_\jph=\frac{(\rho u)^*_\jph}{\rho^*_\jph},\quad\alpha^*_\jph=\left\{\begin{aligned}
&\frac{a^+_\jph}{a^+_\jph-u^*_\jph}&&\mbox{if}~u^*_\jph<0,\\
&\frac{a^-_\jph}{a^-_\jph-u^*_\jph}&&\mbox{otherwise},\end{aligned}\right.\\
&q^\rho_\jph={\rm minmod}\Big(\big(u^*_\jph-a^-_\jph\big)\big(\rho^*_\jph-\rho^-_\jph\big),
\big(a^+_\jph-u^*_\jph\big)\big(\rho^+_\jph-\rho^*_\jph\big)\Big).
\end{aligned}
\label{2.20}
\end{equation}
Note that all of the indexed quantities in \eref{2.18}--\eref{2.20} are time dependent, but from now on we will omit this dependence for the sake of brevity.
\begin{rmk}
As in \cite{KX_22}, the computation of numerical fluxes in \eref{2.18} should be desingularized to avoid division by zero or very small
numbers. If $a^+_\jph<\varepsilon$ and $a^-_\jph>-\varepsilon$ for a small positive $\varepsilon$, we replace the fluxes
$\bm{{\cal F}}_\jph$ with
\begin{equation*}
\bm{{\cal F}}_\jph=\frac{\bm F\big(\bm U^-_\jph\big)+\bm F\big(\bm U^+_\jph\big)}{2}.
\end{equation*}
\end{rmk}

\subsection{New LDCU Scheme for the 2-D Euler Equations of Gas Dynamics}
In this section, we extend the modified semi-discrete LDCU scheme from \S\ref{sec222} to the 2-D Euler equations of gas dynamics, which read
as \eref{1.2} with $\bm U=(\rho,\rho u,\rho v,E)^\top$, $\bm F(\bm U)=(\rho u,\rho u^2+p,\rho uv,u(E+p))^\top$, and
$\bm G(\bm U)=(\rho v,\rho uv,\rho v^2+p,v(E+p))^\top$, where the notations are the same as in the 1-D case except for that now $u$ and $v$
are the $x$- and $y$-velocities, respectively. The system is closed using the EOS for the ideal gas
$\,p=(\gamma-1)\big[E-\frac{\rho}{2}(u^2+v^2)\big]$.

We first introduce a uniform mesh consisting of the finite-volume cells $C_{j,k}:=[x_\jmh,x_\jph]\times[y_\kmh,y_\kph]$ of the uniform size
with $x_\jph-x_\jmh\equiv\dx$ and $y_\kph-y_\kmh\equiv\dy$, $j=1,\ldots,N_x$, $k=1,\ldots,N_y$. We assume that at certain time level $t$, an
approximate solution, realized in terms of the cell averages $\,\xbar{\bm U}_{j,k}$, is available. These cell averages are then evolved in
time by solving the following system of ODEs:
\begin{equation}
\frac{{\rm d}}{{\rm d}t}\,\xbar{\bm U}_{j,k}=-\frac{\bm{{\cal F}}_{\jph,k}-\bm{{\cal F}}_{\jmh,k}}{\dx}-
\frac{\bm{{\cal G}}_{j,\kph}-\bm{{\cal G}}_{j,\kmh}}{\dy},
\label{4.3}
\end{equation}
where the $x$- and $y$-numerical fluxes are
\begin{align}
&\begin{aligned}
\bm{{\cal F}}_{\jph,k}&=\frac{a^+_{\jph,k}\bm F(\bm U^-_{\jph,k})-a^-_{\jph,k}\bm F(\bm U^+_{\jph,k})}{a^+_{\jph,k}-a^-_{\jph,k}}+
\frac{a^+_{\jph,k}a^-_{\jph,k}}{a^+_{\jph,k}-a^-_{\jph,k}}\left(\bm U^+_{\jph,k}-\bm U^-_{\jph,k}\right)\\
&+\bm q_{\jph,k},
\end{aligned}\label{4.4}\\
&\begin{aligned}
\bm{{\cal G}}_{j,\kph}&=\frac{b^+_{j,\kph}\bm G(\bm U^-_{j,\kph})-b^-_{j,\kph}\bm G(\bm U^+_{j,\kph})}{b^+_{j,\kph}-b^-_{j,\kph}}+
\frac{b^+_{j,\kph}b^-_{j,\kph}}{b^+_{j,\kph}-b^-_{j,\kph}}\left(\bm U^+_{j,\kph}-\bm U^-_{j,\kph}\right)\\
&+\bm q_{j,\kph}.
\end{aligned}\label{4.5}
\end{align}

To obtain $\bm U^\pm_{\jph,k}$ and $\bm U^\pm_{j,\kph}$ in \eref{4.4}--\eref{4.5}, we reconstruct the second-order piecewise linear
interpolant $\sum_{j,k}\left[\,\xbar{\bm U}_{j,k}+(\bm U_x)_{j,k}(x-x_j)+(\bm U_y)_{j,k}(y-y_k)\right]{\cal X}_{C_{j,k}}(x,y)$ where
$(\bm U_x)_{j,k}$ and $(\bm U_y)_{j,k}$ are the slopes which are supposed to be computed using a nonlinear limiter to ensure a
non-oscillatory nature of the reconstruction. In the numerical experiments reported in \S\ref{sec3}, we have used the generalized minmod
limiter with the minmod parameter $\theta=1.3$. We then follow \cite[\S3.2]{KX_22} to evaluate the corresponding point values
$u^\pm_{\jph,k}$, $u^\pm_{j,\kph}$, $v^\pm_{\jph,k}$, $v^\pm_{j,\kph}$, $p^\pm_{\jph,k}$, and $p^\pm_{j,\kph}$, and estimate the one-sided
local speeds of propagation $a^\pm_{\jph,k}$ and $b^\pm_{j,\kph}$.

\subsubsection{``Built-in'' Anti-Diffusion}
In this section, we discuss the derivation of the ``built-in'' anti-diffusion terms $\bm q_{\jph,k}$ and $\bm q_{j,\kph}$ in
\eref{4.4}--\eref{4.5} in a ``dimension-by-dimension'' manner following the idea introduced in \cite{KX_22}.

In order to derive the formula for $\bm q_{\jph,k}$, we consider the 1-D restriction of the 2-D system \eref{1.2} along the lines $y=y_k$:
\begin{equation}
\bm U_t(x,y_k,t)+\bm F\big(\bm U(x,y_k,t)\big)_x=\bm0,\quad k=1,\ldots,N_y.
\label{3.17}
\end{equation}
We then go through all of the steps in the derivation of the 1-D fully discrete scheme for the systems in \eref{3.17} following
\cite[\S3.2]{KX_22} and \S\ref{sec2} up to \eref{2.2.6}, which now reads as
\begin{equation*}
\widetilde{\bm U}^{\,\rm int}_{\jph,k}(x,y_k)=\left\{\begin{aligned}
\xbar{\bm U}^{\,\rm int,L}_{\jph,k},\quad x<x_\jph+\widetilde u_{\jph,k}\dt^n,\\
\xbar{\bm U}^{\,\rm int,R}_{\jph,k},\quad x>x_\jph+\widetilde u_{\jph,k}\dt^n,
\end{aligned}\right.
\end{equation*}
and the corresponding local conservation requirements \eref{2.2.7} become
\begin{equation*}
\big(a^+_{\jph,k}-\widetilde u_{\jph,k}\big)\,\xbar{\bm U}^{\,\rm int,R}_{\jph,k}+
\big(\widetilde u_{\jph,k}-a^-_{\jph,k}\big)\,\xbar{\bm U}^{\,\rm int,L}_{\jph,k}=
\big(a^+_{\jph,k}-a^-_{\jph,k}\big)\,\xbar{\bm U}^{\,\rm int}_{\jph,k},
\end{equation*}
where we take $\widetilde u_{\jph,k}=u^{\rm int}_{\jph,k}=(\xbar{\rho u})^{\rm int}_{\jph,k}/\,\xbar\rho^{\,\rm int}_{\jph,k}$. We then
proceed as in \cite[\S3.2]{KX_22}, where we enforce the continuity of $u$ and $p$ across the cell interfaces $x=x_\jph$ by setting
\begin{equation*}
\begin{aligned}
\frac{(\xbar{\rho u})^{\,\rm int,L}_{\jph,k}}{\xbar\rho^{\,\rm int,L}_{\jph,k}}&=
\frac{(\xbar{\rho u})^{\,\rm int,R}_{\jph,k}}{\xbar\rho^{\,\rm int,R}_{\jph,k}},\\
\xbar E^{\,\rm int,L}_{\jph,k}-\frac{\big((\xbar{\rho u})^{\,\rm int,L}_{\jph,k}\big)^2+
\big((\xbar{\rho v})^{\,\rm int,L}_{\jph,k}\big)^2}{2\,\xbar\rho^{\,\rm int,L}_{\jph,k}}
&=\,\xbar E^{\,\rm int,R}_{\jph,k}-\frac{\big((\xbar{\rho u})^{\,\rm int,R}_{\jph,k}\big)^2+
\big((\xbar{\rho v})^{\,\rm int,R}_{\jph,k}\big)^2}{2\,\xbar\rho^{\,\rm int,R}_{\jph,k}},
\end{aligned}
\end{equation*}
and enforce sharp (yet, non-oscillatory) jumps of the $\rho$- and $\rho v$-components. This leads to the following formulae analogous to
\eref{3.14}:
\begin{equation*}
\begin{aligned}
\xbar\rho^{\,\rm int,L}_{\jph,k}&=\,\xbar\rho^{\,\rm int}_{\jph,k}+\frac{\delta^\rho_{\jph,k}}{a^{\rm int,-}_{\jph,k}},\quad&
\xbar\rho^{\,\rm int,R}_{\jph,k}&=\,\xbar\rho^{\,\rm int}_{\jph,k}+\frac{\delta^\rho_{\jph,k}}{a^{\rm int,+}_{\jph,k}},\\[0.5ex]
(\xbar{\rho v})^{\,\rm int,L}_{\jph,k}&=(\xbar{\rho v})^{\,\rm int}_{\jph,k}+\frac{\delta^{\rho v}_{\jph,k}}{a^{\rm int,-}_{\jph,k}},\quad&
(\xbar{\rho v})^{\,\rm int,R}_{\jph,k}&=(\xbar{\rho v})^{\,\rm int}_{\jph,k}+\frac{\delta^{\rho v}_{\jph,k}}{a^{\rm int,+}_{\jph,k}},
\end{aligned}
\end{equation*}
where $a^{\rm int,\pm}_{\jph,k}:=a^\pm_{\jph,k}-u^{\rm int}_{\jph,k}$, and
\begin{equation*}
\begin{aligned}
&\delta^\rho_{\jph,k}={\rm minmod}\left(-a^{\rm int,-}_{\jph,k}
\big[\,\xbar\rho^{\,\rm int}_{\jph,k}-\big(\rho^{\rm int}_{\jph,k}\big)_\ell\big],
a^{\rm int,+}_{\jph,k}\big[\big(\rho^{\rm int}_{\jph,k}\big)_r-\,\xbar\rho^{\,\rm int}_{\jph,k}\big]\right),\\
&\delta^{\rho v}_{\jph,k}={\rm minmod}\left(-a^{\rm int,-}_{\jph,k}
\big[(\xbar{\rho v})^{\,\rm int}_{\jph,k}-\big((\rho v)^{\rm int}_{\jph,k}\big)_\ell\big],
a^{\rm int,+}_{\jph,k}\big[\big((\rho v)^{\rm int}_{\jph,k}\big)_r-(\xbar{\rho v})^{\,\rm int}_{\jph,k}\big]\right),
\end{aligned}
\end{equation*}
and the point values $\big(\rho^{\rm int}_{\jph,k}\big)_\ell$, $\big(\rho^{\rm int}_{\jph,k}\big)_r$,
$\big((\rho v)^{\rm int}_{\jph,k}\big)_\ell$, and $\big((\rho v)^{\rm int}_{\jph,k}\big)_r$ were introduced in \cite[Eq. (2.14)]{KX_22}.

Next, we proceed as in \S\ref{sec22} and complete the derivation of the fully discrete scheme (not shown here for the sake of brevity), and
after this, we pass to the semi-discrete limit and end up with the new LDCU flux \eref{4.4} with the following ``built-in'' anti-diffusion
term:
\begin{equation*}
\bm q_{\jph,k}=\alpha^*_{\jph,k}\big(q^\rho_{\jph,k},u^*_{\jph,k}q^\rho_{\jph,k},q^{\rho v}_{\jph,k},q^E_{\jph,k}\big)^\top.
\end{equation*}
Here,
$$
\begin{aligned}
&\bm U^*_{\jph,k}=\frac{a^+_{\jph,k}\bm U^+_{\jph,k}-a^-_{\jph,k}\bm U^-_{\jph,k}-\big[\bm F(\bm U^+_{\jph,k})-\bm F(\bm U^-_{\jph,k})\big]}
{a^+_{\jph,k}-a^-_{\jph,k}},\quad u^*_{\jph,k}=\frac{(\rho u)^*_{\jph,k}}{\rho^*_{\jph,k}},\\
&q^\rho_{\jph,k}={\rm minmod}\left(-a^{*,-}_{\jph,k}\big(\rho^*_{\jph,k}-\rho^-_{\jph,k}\big),
a^{*,+}_{\jph,k}\big(\rho^+_{\jph,k}-\rho^*_{\jph,k}\big)\right),\\
&q^{\rho v}_{\jph,k}={\rm minmod}\left(-a^{*,-}_{\jph,k}\big((\rho v)^*_{\jph,k}-(\rho v)^-_{\jph,k}\big),
a^{*,+}_{\jph,k}\big((\rho v)^+_{\jph,k}-(\rho v)^*_{\jph,k}\big)\right),\\
&q^E_{\jph,k}=\frac{a^{*,+}_{\jph,k}a^{*,-}_{\jph,k}}{a^+_{\jph,k}-a^-_{\jph,k}}
\left\{\frac{\Bigg((\rho v)^*_{\jph,k}+\dfrac{q^{\rho v}_{\jph,k}}{a^{*,+}_{\jph,k}}\Bigg)^2}
{2\Bigg(\rho^*_{\jph,k}+\dfrac{q^\rho_{\jph,k}}{a^{*,+}_{\jph,k}}\Bigg)}-
\frac{\Bigg((\rho v)^*_{\jph,k}+\dfrac{q^{\rho v}_{\jph,k}}{a^{*,-}_{\jph,k}}\Bigg)^2}
{2\Bigg(\rho^*_{\jph,k}+\dfrac{q^\rho_{\jph,k}}{a^{*,-}_{\jph,k}}\Bigg)}\right\}+\frac{\big(u^*_{\jph,k}\big)^2}{2}q^\rho_{\jph,k},
\end{aligned}
$$
and
\begin{equation*}
\alpha^*_{\jph,k}=\left\{\begin{aligned}&\frac{a^+_{\jph,k}}{a^{*,+}_{\jph,k}}&&\mbox{if}~u^*_{\jph,k}<0,\\
&\frac{a^-_{\jph,k}}{a^{*,-}_{\jph,k}}&&\mbox{otherwise},\end{aligned}\right.\qquad a^{*,\pm}_{\jph,k}=a^\pm_{\jph,k}-u^*_{\jph,k}.
\end{equation*}
Similarly, the ``built-in'' anti-diffusion term in the new $y$-directional LDCU flux \eref{4.5} is
\begin{equation*}
\bm q_{j,\kph}=\alpha^*_{j,\kph}\big(q^\rho_{j,\kph},q^{\rho u}_{j,\kph},v^*_{j,\kph}q^\rho_{j,\kph},q^E_{j,\kph}\big)^\top
\end{equation*}
with
$$
\begin{aligned}
&\bm U^*_{j,\kph}=\frac{b^+_{j,\kph}\bm U^+_{j,\kph}-b^-_{j,\kph}\bm U^-_{j,\kph}-\big[\bm G(\bm U^+_{j,\kph})-\bm G(\bm U^-_{j,\kph})\big]}
{b^+_{j,\kph}-b^-_{j,\kph}},\quad v^*_{j,\kph}=\frac{(\rho v)^*_{j,\kph}}{\rho^*_{j,\kph}},\\
&q^\rho_{j,\kph}={\rm minmod}\left(-b^{*,-}_{j,\kph}\big(\rho^*_{j,\kph}-\rho^-_{j,\kph}\big),
b^{*,+}_{j,\kph}\big(\rho^+_{j,\kph}-\rho^*_{j,\kph}\big)\right),\\
&q^{\rho u}_{j,\kph}={\rm minmod}\left(-b^{*,-}_{j,\kph}\big((\rho u)^*_{j,\kph}-(\rho u)^-_{j,\kph}\big),
b^{*,+}_{j,\kph}\big((\rho u)^+_{j,\kph}-(\rho u)^*_{j,\kph}\big)\right),\\
&q^E_{j,\kph}=\frac{b^{*,+}_{j,\kph}b^{*,-}_{j,\kph}}{b^+_{j,\kph}-b^-_{j,\kph}}
\left\{\frac{\Bigg((\rho u)^*_{j,\kph}+\dfrac{q^{\rho u}_{j,\kph}}{b^{*,+}_{j,\kph}}\Bigg)^2}
{2\Bigg(\rho^*_{j,\kph}+\dfrac{q^\rho_{j,\kph}}{b^{*,+}_{j,\kph}}\Bigg)}-
\frac{\Bigg((\rho u)^*_{j,\kph}+\dfrac{q^{\rho u}_{j,\kph}}{b^{*,-}_{j,\kph}}\Bigg)^2}
{2\Bigg(\rho^*_{j,\kph}+\dfrac{q^\rho_{j,\kph}}{b^{*,-}_{j,\kph}}\Bigg)}\right\}+\frac{\big(v^*_{j,\kph}\big)^2}{2}q^\rho_{j,\kph},
\end{aligned}
$$
and
\begin{equation*}
\alpha^*_{j,\kph}=\left\{\begin{aligned}&\frac{b^+_{j,\kph}}{b^{*,+}_{j,\kph}}&&\mbox{if}~v^*_{j,\kph}<0,\\
&\frac{b^-_{j,\kph}}{b^{*,-}_{j,\kph}}&&\mbox{otherwise},\end{aligned}\right.\qquad b^{*,\pm}_{j,\kph}=b^\pm_{j,\kph}-v^*_{j,\kph}.
\end{equation*}
\begin{rmk}
As in the 1-D case, the computation of numerical fluxes in \eref{4.4} and \eref{4.5} should be desingularized to avoid division by zero or
very small numbers:

\noindent
$\bullet$ If $a^+_{\jph,k}<\varepsilon$ and $a^-_{\jph,k}>-\varepsilon$ for a small positive $\varepsilon$, we replace the flux
$\bm{{\cal F}}_{\jph,k}$ with
\begin{equation*}
\bm{{\cal F}}_{\jph,k}=\frac{\bm F\big(\bm U^-_{\jph,k}\big)+\bm F\big(\bm U^+_{\jph,k}\big)}{2};
\end{equation*}

\noindent
$\bullet$ If $b^+_{j,\kph}<\varepsilon$ and $b^-_{j,\kph}>-\varepsilon$, we replace the flux $\bm{{\cal G}}_{j,\kph}$ with
\begin{equation*}
\bm{{\cal G}}_{j,\kph}=\frac{\bm G\big(\bm U^-_{j,\kph}\big)+\bm G\big(\bm U^+_{j,\kph}\big)}{2}.
\end{equation*}
\end{rmk}

\section{Numerical Examples}\label{sec3}
In this section, we apply the developed schemes to several initial-boundary value problems for the 1-D and 2-D Euler equations of gas
dynamics (with $\gamma=1.4$) and compare the performance of the new and original second-order LDCU schemes, which will be referred to as the
NEW and OLD schemes.

For time integration of the ODE systems \eref{2.17} and \eref{4.3}, we have used the three-stage third-order strong stability preserving
(SSP) Runge-Kutta method (see, e.g.,\cite{Gottlieb11,Gottlieb12}) with the CFL number 0.475. We have taken the small desingularization
parameter $\varepsilon=10^{-12}$.

\subsubsection*{Example 1---Shock-Entropy Problem}
In the first example taken from \cite{Shu88}, we consider the shock-entropy wave interaction problem with the  following initial condition:
\begin{equation*}
\big(\rho(x,0),u(x,0),p(x,0)\big)=\begin{cases}(1.51695,0.523346,1.805),&x<-4.5,\\(1+0.1\sin(20x),0,1),&x>-4.5,\end{cases}
\end{equation*}
which corresponds to a forward-facing shock wave of Mach number 1.1 interacting with high-frequency density perturbations, that is, as the
shock wave moves, the perturbations spread ahead.

\begin{figure}[ht!]
\centerline{\includegraphics[trim=0.7cm 0.4cm 1.2cm 0.6cm, clip, width=6.5cm]{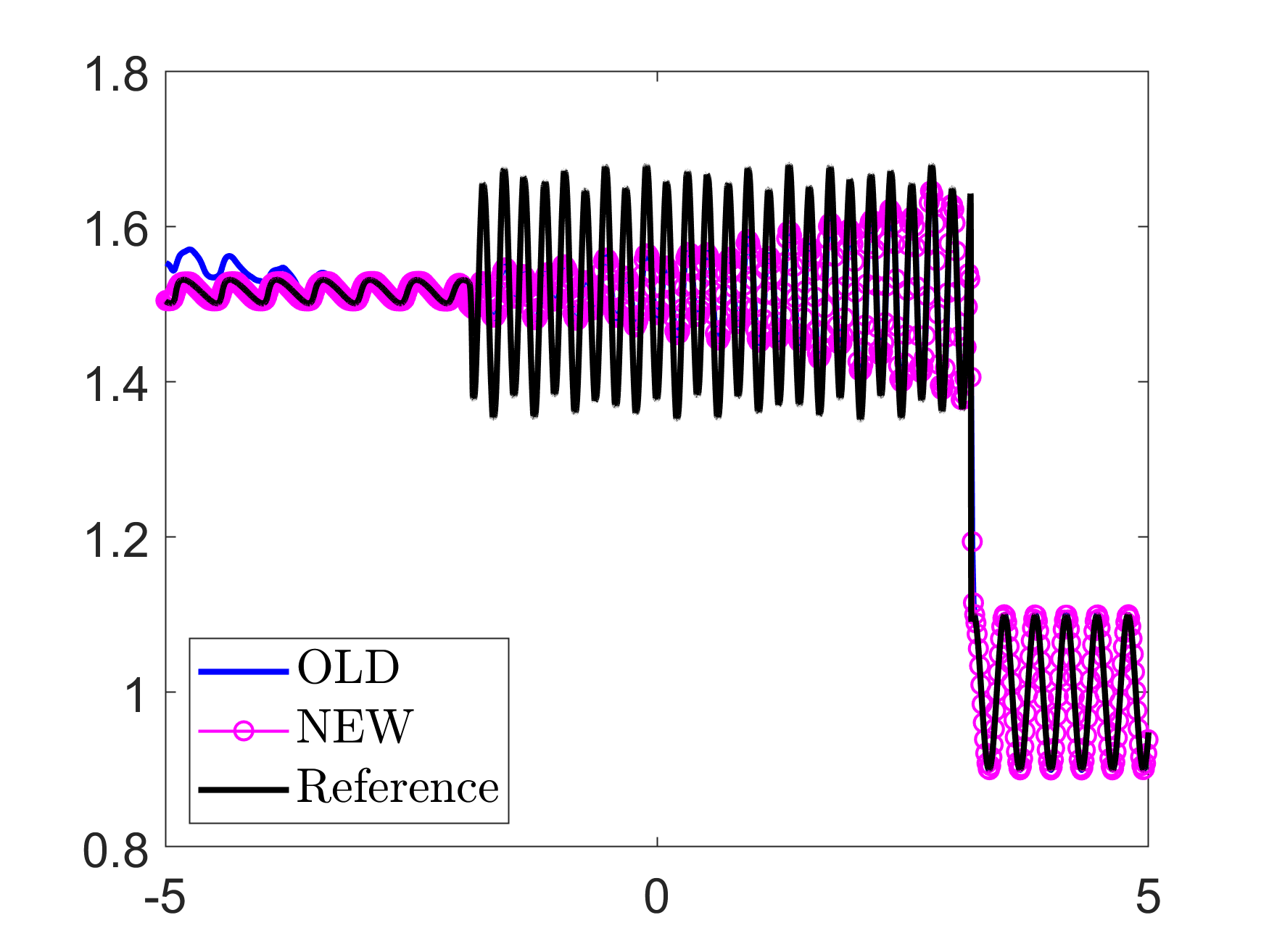}\hspace*{1.0cm}
            \includegraphics[trim=0.7cm 0.4cm 1.2cm 0.6cm, clip, width=6.5cm]{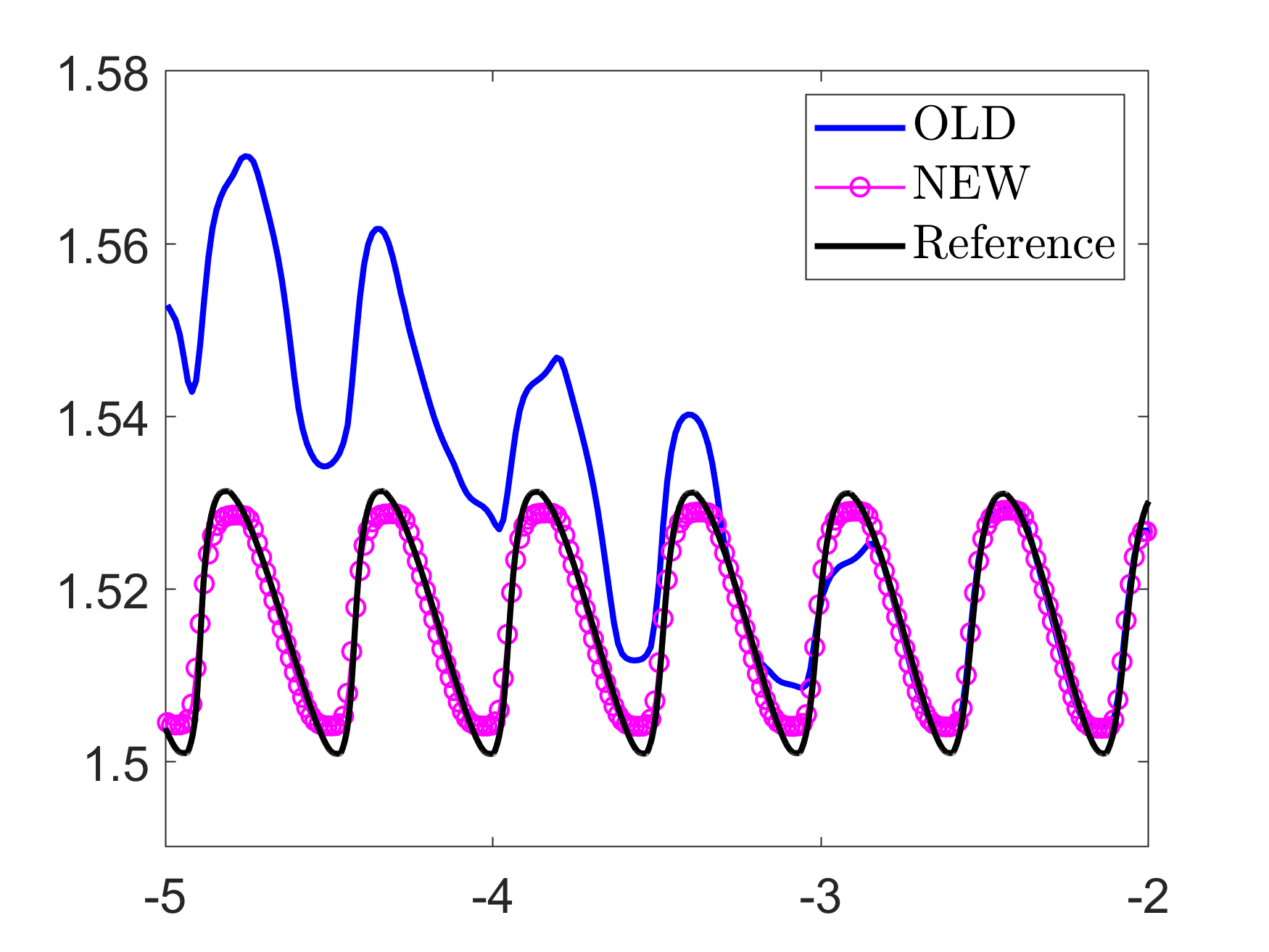}}
\caption{\sf Example 1: Density $\rho$ computed by the OLD and NEW schemes (left) and zoom at $x\in[-5,-2]$ (right).\label{fig1}}
\end{figure}

We compute the numerical solution using both the NEW and OLD schemes in the computational domain $[-5,5]$ on a uniform mesh with $\dx=1/80$.
We impose the free boundary conditions by simply setting $\,\xbar{\bm U}_0:=\,\xbar{\bm U}_1$ and
$\,\xbar{\bm U}_{N+1}:=\,\xbar{\bm U}_N$ in the ghost cells $C_0$ and $C_{N+1}$ on the left and on the right, respectively. The numerical
results at time $t=5$ are presented in Figure \ref{fig1} along with the corresponding reference solution computed by the NEW scheme on a
much finer mesh with $\dx=1/800$. As one can see, the numerical solution computed by the OLD scheme is very inaccurate near the left
boundary of the computational domain, while the NEW scheme accurately captures the solution throughout the entire computational domain.

\subsubsection*{Example 2---Stationary Contact Wave, Traveling Shock, and Rarefaction Wave}
In the second example taken from \cite{Kurganov07}, the initial conditions,
\begin{equation*}
\big(\rho(x,0),u(x,0),p(x,0)\big)=\begin{cases}(1,-19.59745,1000)&\mbox{if}~x<0.8,\\(1,-19.59745,0.01)&\mbox{otherwise},\end{cases}
\end{equation*}
are prescribed in the computational domain $[0,1]$ subject to the free boundary conditions.

\begin{figure}[ht!]
\centerline{\includegraphics[trim=1.4cm 0.4cm 0.8cm 0.4cm, clip, width=6.5cm]{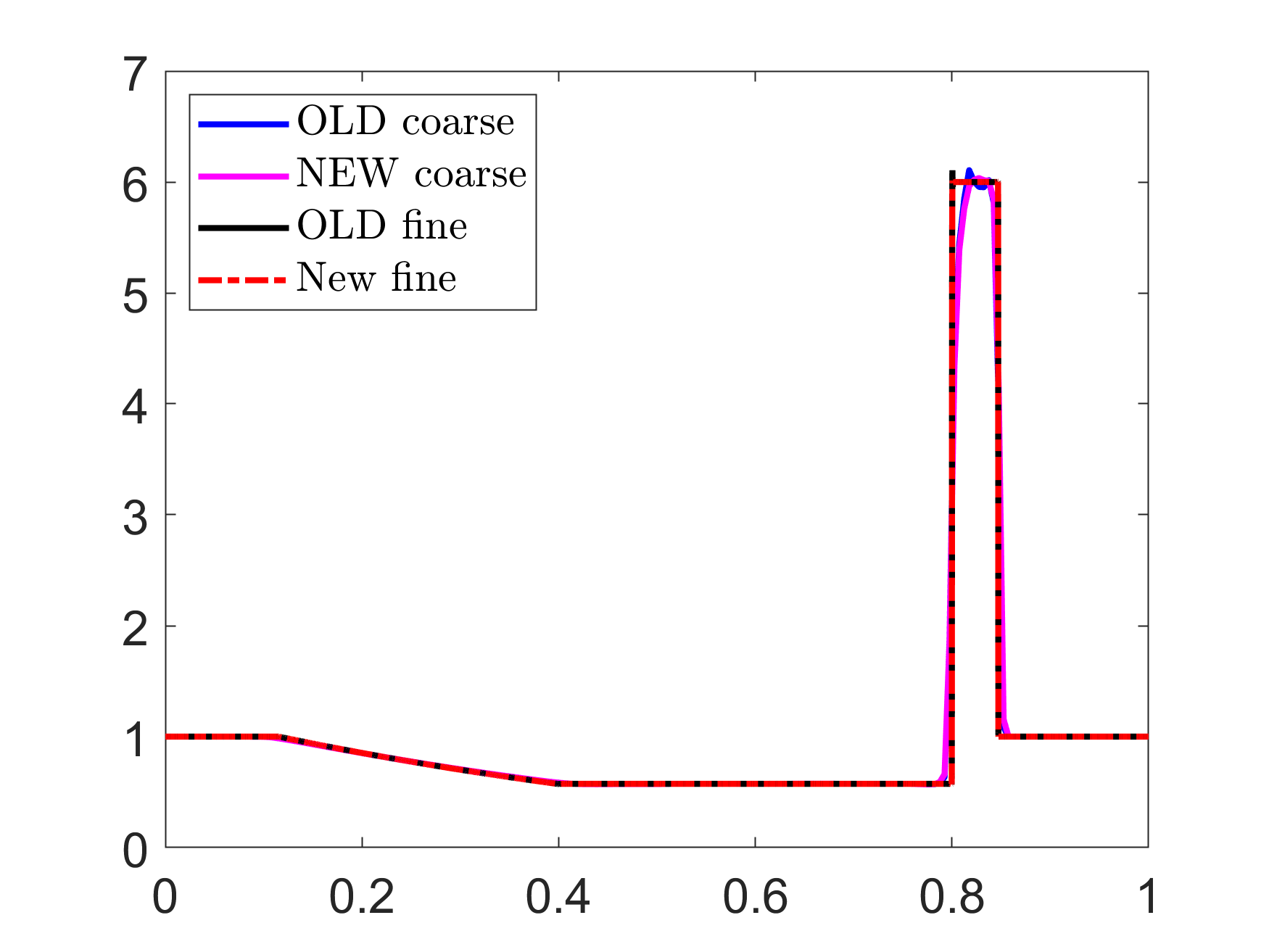}\hspace*{1.0cm}
            \includegraphics[trim=1.4cm 0.4cm 0.8cm 0.4cm, clip, width=6.5cm]{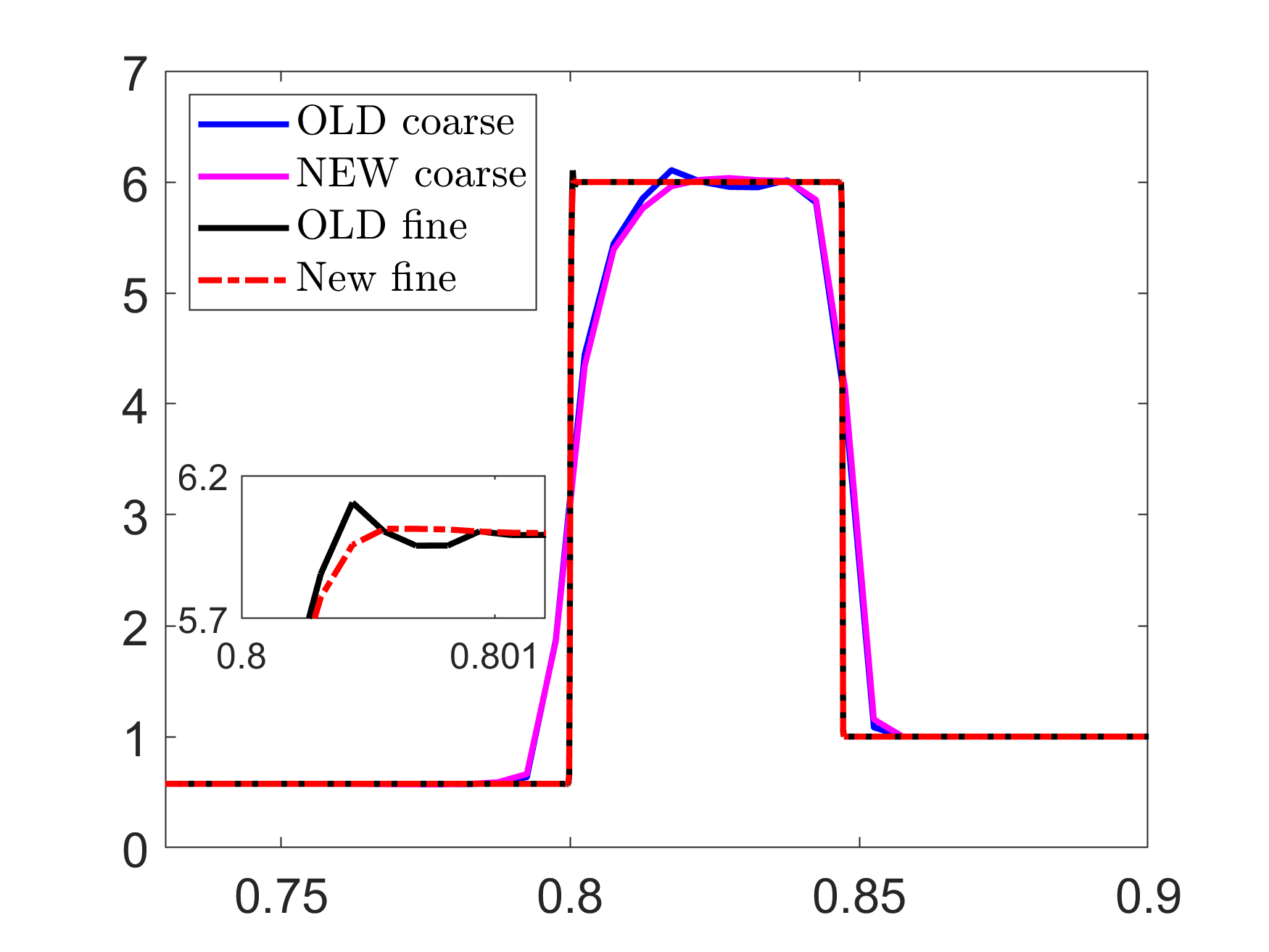}}
\caption{\sf Example 2: Density $\rho$ computed by the OLD and NEW schemes on two uniform meshes (left) and zoom at $x\in[0.7,0.9]$ and $x\in[0.8,0.801]$ (right).\label{fig2}}
\end{figure}

We compute the numerical solutions until the final time $t=0.012$ using both the NEW and OLD schemes on two uniform meshes, the coarse and
fine ones with  $\dx=1/200$ and $1/8000$, respectively. The numerical results, plotted in Figure \ref{fig2}, show that the NEW and OLD
schemes achieve similar resolutions of both shock and contact waves. At the same time, the numerical results computed by the NEW scheme is
non-oscillatory, while the OLD scheme solutions contain small oscillations, whose magnitude does not decay when the mesh is refined.

\subsubsection*{Example 3---2-D Riemann Problem} In the first 2-D example, we consider Configuration 3 of the 2-D Riemann problems taken
from \cite{KT2002}; see also \cite{Z2001,S1993,SCG1993}. The initial conditions,
\begin{equation*}
(\rho(x,y,0),u(x,y,0),v(x,y,0),p(x,y,0))=\left\{\begin{aligned}
&(1.5,0,0,1.5),&&x>1,~y>1,\\
&(0.5323,1.206,0,0.3),&&x<1,~y>1,\\
&(0.138,1.206,1.206,0.029),&&x<1,~y<1,\\
&(0.5323,0,1.206,0.3),&&x>1,~y<1,
\end{aligned}\right.
\end{equation*}
are prescribed in the computational domain $[0,1.2]\times[0,1.2]$ subject to the free boundary conditions.

\begin{figure}[ht!]
\centerline{\includegraphics[trim=5.4cm 3.7cm 2.4cm 2.3cm,clip,width=15cm]{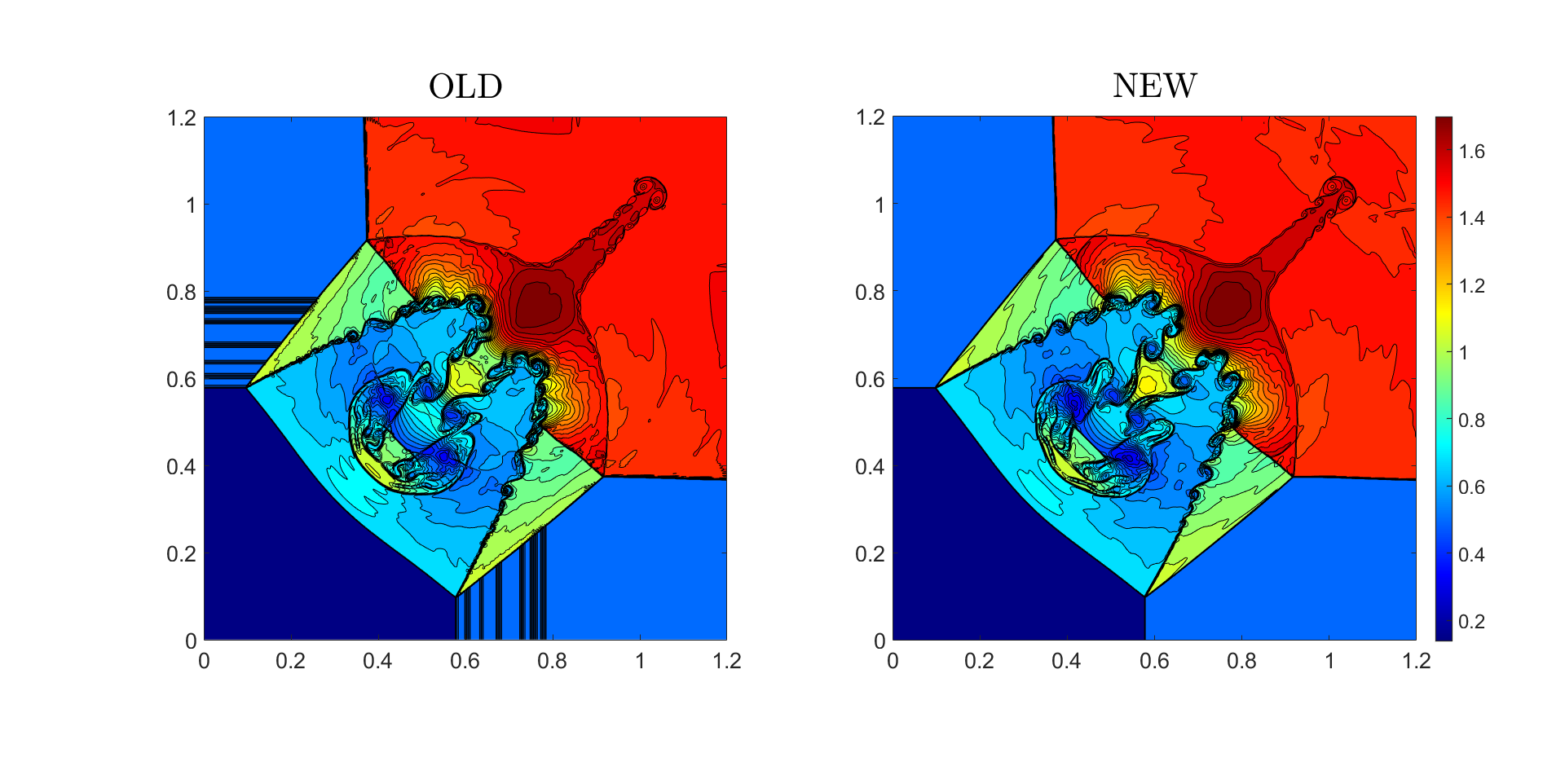}}
\caption{\sf Example 3: Density ($\rho$) computed by the OLD (left) and NEW (right) schemes.\label{fig2.2}}
\end{figure}

We compute the solution using both the NEW and OLD schemes on a uniform mesh with $\dx=\dy=0.001$ until the final time $t=1$ and present the
obtained densities in Figure \ref{fig2.2}. As one can see, the NEW solution is substantially less oscillatory than the OLD one, and the
resolution achieved by the NEW scheme seems to be comparable and even a little better in some parts of the computed solutions.

\subsubsection*{Example 4---Explosion Problem}
In this example, we consider the explosion problem taken from \cite{LW2003} with the following initial conditions:
\begin{equation*}
(\rho(x,y,0),u(x,y,0),v(x,y,0),p(x,y,0))=\left\{\begin{aligned}
&(1,0,0,1),&&x^2+y^2<0.16,\\
&(0.125,0,0,0.1),&&\mbox{otherwise},
\end{aligned}\right.
\end{equation*}
which are prescribed in the computational domain $[0,1.5]\times[0,1.5]$ subject to the solid wall boundary conditions at $x=0$ and $y=0$ and
free boundary conditions at $x=1.5$ and $y=1.5$.

\begin{figure}[ht!]
\centerline{\includegraphics[trim=5.4cm 3.7cm 2.2cm 2.3cm,clip,width=15cm]{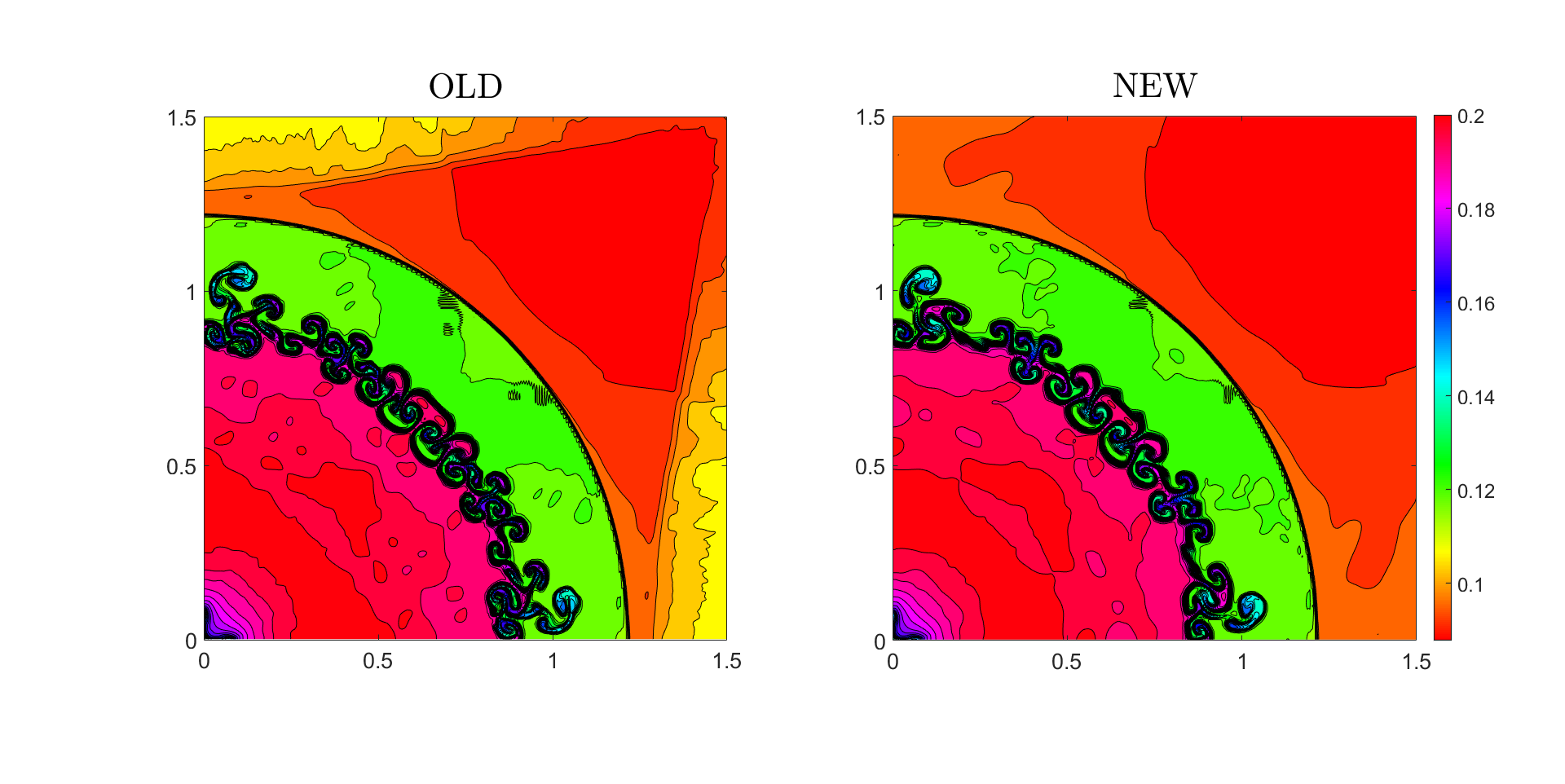}}
\caption{\sf Example 4: Density ($\rho$) computed by the OLD (left) and NEW (right) schemes.\label{fig2.3}}
\end{figure}

We compute the solution using the NEW and OLD schemes on a uniform mesh with $\dx=\dy=3/800$ until the final time $t=3.2$. The obtained
densities are presented in Figure \ref{fig2.3}, where one can clearly see that there are obvious oscillations along the boundaries $x=1.5$
and $y=1.5$ in the numerical results computed by the OLD scheme, while the oscillations are substantially smaller in the numerical results
computed by the NEW scheme. At the same time, in this example, the OLD scheme achieves a slightly better resolution.

\paragraph{Example 5---Implosion Problem}
In the last example, we test the implosion problem also taken from \cite{LW2003}. The initial conditions,
\begin{equation*}
(\rho(x,y,0),u(x,y,0),v(x,y,0),p(x,y,0))=\left\{\begin{aligned}
&(0.125,0,0,0.14),&&\vert x\vert+\vert y\vert<0.15,\\
&(1,0,0,1),&&\mbox{otherwise},
\end{aligned}\right.
\end{equation*}
are prescribed in the computational domain $[0,0.3]\times[0,0.3]$ subject to the solid wall boundary conditions.

\begin{figure}[ht!]
\centerline{\includegraphics[trim=5.4cm 3.7cm 2.4cm 2.3cm,clip,width=15cm]{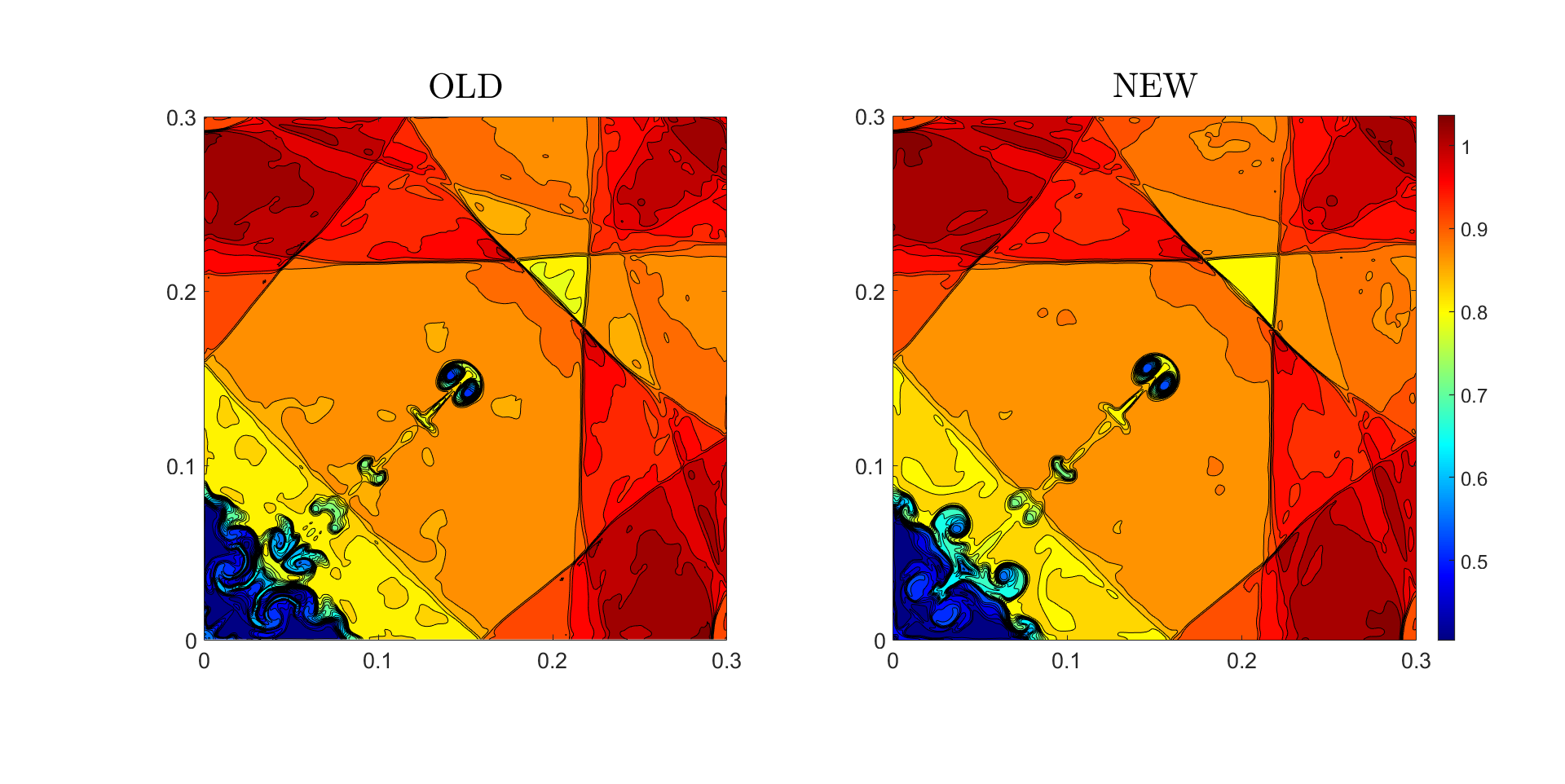}}
\caption{\sf Example 5: Density ($\rho$) computed by the OLD (left) and NEW (right) schemes.\label{fig2.4}}
\end{figure}

We compute the solution using the NEW and OLD schemes on a uniform mesh with $\dx=\dy=1/2000$ until the final time $t=2.5$. The obtained
densities are presented in Figure \ref{fig2.4}. As one can see, the NEW solution is slightly less oscillatory than the OLD one. At the same
time, the jet generated by the NEW scheme propagates to a larger extent than the jet produced by the OLD scheme, which demonstrates that in
this example, the NEW scheme achieves slightly higher resolution than the OLD scheme.

\begin{DA}
\paragraph{Funding.} The work of A. Kurganov was supported in part by NSFC grant 12171226, and by the fund of the Guangdong Provincial Key Laboratory of Computational Science and Material Design (No. 2019B030301001).

\paragraph{Conflicts of interest.} On behalf of all authors, the corresponding author states that there is no conflict of interest.

\paragraph{Data and software availability.} The data that support the findings of this study and FORTRAN codes developed by the authors and used to obtain all of the presented numerical results are available from the corresponding author upon reasonable request.
\end{DA}

\bibliography{reference}
\bibliographystyle{siamnodash}
\end{document}